\documentclass[11pt]{article}
\usepackage{amsfonts,epsf,amsmath,amssymb,tikz}

\newtheorem{theorem}{\bf Theorem}

\newtheorem{lemma}[theorem]{\bf Lemma}

\newcommand{\proof}{\noindent{\bf Proof.\ }}
\newcommand{\qed}{\hfill $\square$ \bigskip}
\newcommand{\bnd}{\frac{n}{g}(m-n+1)}

\author{
Jernej Azarija \\
Institute of Mathematics, Physics and Mechanics \\
Jadranska 19, 1000 Ljubljana, Slovenia \\
jernej.azarija@gmail.com
\and
Sandi Klav\v zar \\
Faculty of Mathematics and Physics, University of Ljubljana \\
Jadranska 19, 1000 Ljubljana, Slovenia \\
and \\
Faculty of Natural Sciences and Mathematics, University of Maribor \\
Koro\v ska 160, 2000 Maribor, Slovenia \\
sandi.klavzar@fmf.uni-lj.si 
}

\date{\today}

\title{Moore graphs and cycles are extremal graphs for convex cycles}

\begin{document}
\maketitle

\begin{abstract}
Let $\rho(G)$ denote the number of convex cycles of a simple graph $G$ of order $n$, size $m$, and girth $3 \leq g \leq n$. It is proved that $\rho(G) \leq \bnd $ and that equality holds if and only if $G$ is an even cycle or a Moore graph. The equality also holds for a possible Moore graph of diameter 2 and degree 57 thus giving a new characterization of Moore graphs.
\end{abstract}


\noindent {\bf AMS Subj. Class. (2010)}: 05C75, 05C12

\section{Introduction}

Convexity is a central notion in the theory of discrete metric spaces~\cite{Vel-1993}.  
In graph theory, convex subgraphs and in particular convex cycles are often employed
to unveil additional structure of the studied graphs. Recall that a subgraph $H$ of a graph $G$ is {\em convex} if for any $u,v\in V(H)$, every shortest $u,v$-path in $G$ lies completely in $H$. In particular, if $H$ is convex in $G$, then $d_H(u,v) = d_G(u,v)$ holds for any $u,v\in V(H)$, where $d_G$ denotes the usual shortest path distance in $G$.  

Convex subgraphs are indispensable in the study of (Cartesian) graph products. 
Extending a result of Vanden Cruyce~\cite{Vanden-1982} for hypercubes, Egawa~\cite{Egawa-1986}
characterized Cartesian products of complete graphs by convex subgraphs. 
Similarly, Chepoi~\cite{Chepoi-1988} characterized isometric subgraphs of 
Cartesian products of complete graphs via convexity of certain 
subgraphs. In~\cite{Bandelt-1996a} weak Cartesian products of trees are characterized 
among median graphs by the property that $K_{2,3}$ minus an edge is not a convex subgraph. 
For additional aspects on the convexity of graph products see~\cite{Hammack-2011}. 
For instance, the book contains a short proof of the classical 
unique prime factorization theorem with respect to the Cartesian product,
where convexity is the key tool for the short proof.  

Convex subgraphs are even more important in understanding the structure of isometric 
subgraphs of hypercubes, graphs known as partial cubes. (Recall that 
median graphs form a distinguished subclass of partial cubes.) It all started
with the seminal paper of Djokovi\v c~\cite{Djokovic-1973} in which partial
cubes are characterized among bipartite graph with the convexity of subgraphs
induced by vertices closer to one endpoint of an edge than to the other. 
Later, Bandelt and Chepoi~\cite{Bandelt-1996b} characterized acyclic cubical complexes
among median graphs by forbidden convex subgraphs. These graphs were further  
characterized as the graphs for which $-2$ is a zero of the cube polynomial of 
an arbitrary 2-connected convex subgraph~\cite{Bresar-2006}. 

Among convex subgraphs, convex cycles are frequently studied. 
In~\cite{Polat-2009a} Polat proved that a netlike partial cube is prism-retractable 
if and only if it contains at most one convex cycle of length greater than 4 while 
in~\cite{Polat-2009b} he showed that any netlike partial cube that is without an isometric 
ray contains a convex cycle or a finite hypercube which is fixed by every automorphism. 
Parallel to the first mentioned Polat's result it was proved in~\cite{Klavzar-2012b} that 
a partial cube is almost-median if and only if it contains no convex cycle of length 
greater than 4. Very recently the convex excess of a graph was 
introduced as the sum of contributions of all of its convex cycles and used to obtain an 
inequality involving the order, the size, the isometric dimension, and the convex excess 
of an arbitrary partial cube~\cite{Klavzar-2012a}. 

Here we consider convex cycles from an extremal point of view: what is the largest
number of convex cycles that a given graph can have? We became interested in this 
question because of the recent paper~\cite{Hellmuth-2011} by Hellmuth, Leydold, and
Stadler in which convex cycle bases are studied. Along the way they also proved that 
a graph $G$ of order $n$ and size $m$ contains at most $nm$ convex cycles. In this paper
we strengthen this by proving the following result. 

\begin{theorem}
\label{thm:main}
Let $G$ be a simple graph of order $n$, size $m$, and girth $g\ge 3$. Then $G$ contains at most 
$$\bnd$$
convex cycles. Moreover, equality holds if and only if $G$ is an even cycle or a Moore graph.
\end{theorem}

Recall that a {\em Moore graph} is a graph with the maximum possible number of 
vertices that a given graph with prescribed
maximum degree and diameter can have. Equivalently, a Moore graph can also be 
defined as a graph with diameter $r$ and girth $2r+1$, cf.~\cite[p.~90]{Godsil}.  
Singleton~\cite{Singleton} proved that Moore graphs are regular, 
see~\cite[Lemma~5.8.1]{Godsil} for an elegant proof.  
The only Moore graphs that exist are complete graphs, 
odd cycles, the Petersen graph, the Hoffman-Singleton graph, and possibly a Moore
graph of diameter 2 and degree 57~\cite{Hoffman-1960,bannai-1973,Damerell-1973}.
The existence of a latter Moore graph is a big open problem. 
As a by-product of teh proof of Theorem~\ref{thm:main}, we also get: 

\begin{theorem}
\label{thm:main2}
Let $G$ be a simple graph of order $n$, size $m$, and girth $g=2r+1$. Then $G$ is a Moore graph if and only if  the number of $(2r+1)$-cycles in $G$ is $\frac{n}{2r+1}(m-n+1)$.
\end{theorem}
  
For detailed information on Moore graphs and related classes of graphs see the survey~\cite{Miller-2005}. 
The recent paper~\cite{Macaj-2010} contains further insights into a missing Moore graph.
In particular it is proved that the order of the automorphism group of such a graph is 
at most 375 thus significantly extending the fact that it is not vertex-transitive as
proved Graham Higman in a series of lectures, cf.~\cite[Theorem 3.13]{Cameron}. 
On the other hand \v Siagiov\'{a} and \v Sir\'{a}\v{n}~\cite{Siagiova-2012} proved that for an infinite set of degrees $r$ there exist vertex-transitive graphs of degree $r$, diameter 2, 
and order close to the Moore bound. 

The next section contains the proof of Theorems~\ref{thm:main} and~\ref{thm:main2}, a concluding remark is given in the final section.

\section{Proof of Theorem~\ref{thm:main}}

This section is organized as follows. We first characterize convex cycles in a way suitable to us. In the following subsection the number of odd convex cycles of a given graph is bounded and proved that precisely the Moore graphs are extremal graphs. In Subsection~\ref{SEC:even} we then prove a corresponding upper bound for even convex cycles while in the last subsection a combined inequality is derived.     

In what follows $G$ will denote a simple graph on $n$ vertices, with $m$ edges, and of girth $g\ge 3$. The following characterization of convex cycles is a modification of a related result proved in~\cite{Hellmuth-2011}.
More precisely, the first part (for odd cycles) is the same, while the second part is modified to serve our purposes. 

\begin{lemma} \label{LM:ccchar}
Let $C$ be a cycle of $G$.  If $|C| = 2k+1$, $k\ge 1$, then $C$ is convex if and only if for every edge $e = xy$ of $C$ there exists a vertex $v \in C$ such that
    \begin{description}
      \item[(i)] $d_G(x,v) = d_G(y,v) = k$, and 
      \item[(ii)] the $x,v$-path (resp. $y,v$-path) on $C$ of length $k$ is a unique shortest $x,v$-path (resp. $y,v$-path) in $G$. 
\end{description}
If $|C| =  2k$, $k\ge 2$, then $C$ is convex if and only if for every vertex $u\in C$ there exists a vertex $v\in C$ such that 
    \begin{description}
        \item[(iii)] $d_G(u,v) = k$,
        \item[(iv)] there are precisely two $u,v$-paths in $G$ of length $k$. 
\end{description}
\end{lemma}

\proof
As mentioned above, we only need to prove the even case. Hence let $|C| = 2k$, $k\ge 2$. It is clear that the two conditions are necessary. Suppose now that for every vertex $u\in C$ there exists a vertex $v\in C$ such that (iii) and (iv) hold. By way of contradiction assume that there are vertices $x,y\in C$ such that there is shortest $x,y$-path $P$ that is not completely contained in $C$. Let $x'$ be the vertex on $C$  with $d_G(x,x')=k$.  By (iv) there are precisely two $x,x'$-paths in $G$ of length $k$ and they are both contained in $C$. Then $y$ belongs to one of these paths, denote it with $Q$. If $P$ is shorter than the length of the $x,y$-subpath of $Q$, then $d_G(x,x') < k$, a contradiction. And if $P$ is of the same length as the $x,y$-subpath of $Q$, then we would have at least three $x,x'$-paths of length $k$, which contradicts (iv) for $x$ and $x'$. 
\qed

For later use we note here that if follows from the first part of Lemma~\ref{LM:ccchar} that in a graph of girth $g=2r+1$ all of its $g$-cycles are convex. 

We will call a pair $(e,v) \in E(G) \times V(G)$ that satisfies conditions (i) and (ii) of Lemma \ref{LM:ccchar} an {\em odd antipodal pair}. Likewise if $(u,v) \in V(G) \times V(G)$ satisfies conditions (iii) and (iv) then we will say that $(u,v)$ is an {\em even antipodal pair}. In cases where the context is clear we will simply say that a pair $(a,b)$ is {\em antipodal} if it is an even or odd antipodal pair.  

Observe that Lemma \ref{LM:ccchar} readily implies that the number of odd convex cycles is $O(nm)$ while the number of even convex cycles is $O(n^2).$ In what follows we give sharper estimates for these two quantities by bounding the number of antipodal pairs.

\subsection{Odd convex cycles}\label{SEC:odd}

\begin{lemma} \label{LM:oddantipodal}
For any vertex $v \in V(G)$ there are at most $m-n+1$ edges $e$ such that $(e,v)$ is an odd antipodal pair.
\end{lemma}

\proof
Let $T$ be a BFS tree of $G$ with root $v$. Then the assertion readily follows from the fact that if $e\in E(T)$, then one endpoint of $e$ is closer to $v$ than the other. Consequently, $(e,v)$ is not antipodal.  
\qed

From Lemma~\ref{LM:oddantipodal} we get an estimate on the number of odd convex cycles in $G$, which we denote by $\rho_o(G)$. 

\begin{lemma}\label{CL:boundodd}
$\rho_o(G) \leq \bnd$.
\end{lemma}

\proof
Suppose that $G$ contains $k$ odd convex cycles. Every convex cycle $C$ determines precisely $|C| \geq g$ antipodal pairs. We select one and assign it to $C$. Doing it for every convex cycle, there are at least $k(g-1)$ antipodal pairs that are not assigned to convex cycles. In addition, by Lemma~\ref{LM:oddantipodal}, a vertex of $G$ does not form an antipodal pair with at least $n-1$ edges. Therefore we have at least $n(n-1)$ non-antipodal pairs. If follows that 
$$ k\le nm - k(g-1) - n(n-1) $$ 
and thus 
$$k \leq \bnd$$ 
as claimed.
\qed

If $G$ is a cycle, then $m = n = g$, thus the bound of Lemma~\ref{CL:boundodd} is sharp for all odd cycles. The same holds for complete graphs $K_n$, $n \geq 3$. Indeed, for $K_n$ we have $g=3$, $m={n \choose 2}$, and any triple of vertices  induces a triangle, hence the assertion follows because $\frac{n}{3}\left( {n \choose 2} - n +1\right) = {n \choose 3}$. We next show that equality in Lemma~\ref{CL:boundodd} holds precisely for the Moore graphs. 

\begin{lemma}\label{CL:charoddbound}
$\rho_o(G) = \bnd$ if and only if $G$ is a Moore graph.
\end{lemma}

\proof
Suppose first that $G$ is a graph that satisfies the equality. Then it follows from Lemma~\ref{CL:boundodd} and its proof that the girth $g$ of $G$ is odd and that all convex cycles of $G$ are of length $g = 2r+1$. Recall from Lemma~\ref{LM:oddantipodal} that a vertex $v \in V(G)$ lies in at most $m-n+1$ antipodal pairs. Since the equality is satisfied for $G$, it follows that every edge which is not on a BFS tree with a root $v$ constitutes  an antipodal pair with $v$. 
In other words every such edge joins two vertices $x,y$ such that $d_G(v,x) = d_G(v,y) = r$. 
Consider now a BFS tree $T$ rooted at $v$ and let $v'$ be a leaf of $T$. Observe that $v'$ has degree at least two in $G$ because $\rho(G) = \rho(G-u)$ holds for any pendant vertex $u$. Hence there is an edge $e$ not in $T$ that is adjacent to $v'$ in $G$. From the above remark it follows that $(e,v)$ is an antipodal pair and therefore $d_G(v,v') = r$. This in turn implies that $G$ has diameter $r$. Since the  girth of $G$ is $2r+1$ we conclude that $G$ is a Moore graph. 

To prove the converse we need to show that every Moore graph satisfies the equality. As already observed, this is the case with odd cycles and complete graphs of order $n \geq 3$. The Petersen graph has girth 5, hence all of its twelve 5-cycles are convex. Since $\frac{10}{5}(15-10+1) = 12$, the bound for the Petersen graph is established. 

It thus remains to show that the Hoffman-Singleton graph $H$ and a possible Moore graph $X$ of diameter 2 and degree 57 also have the claimed property. To show this we use an implication of a result of Harary \cite{Harary} which can be formulated as follows, cf.~\cite[p.~45]{Biggs}. Let $p_G(x)$ be the characteristic polynomial of a graph $G$ of order $n$ and odd girth $g$. Then the number of $g$-cycles of $G$ equals $-c/2$, where $c$ is the coefficient at $x^{n-g}$ in $p_G(x)$. 

The Hoffman-Singleton graph $H$ has 50 vertices, 175 edges, and $p_H(x) = (x-7)(x-2)^{28}(x+3)^{21}$, cf.~\cite{Rowlinson}. Since it has girth $5$ and 
$$\frac{\left( \frac{d^{45}}{dx^{45}} p_H(x) \right)(0)}{45!} = -2520\,,$$
it follows that the number of 5-cycles of $H$ is 1260. Hence the bound of Lemma~\ref{CL:charoddbound} is sharp for $H$.

For the possible Moore graph $X$ it is known that $p_X(x) = (x-57)(x+8)^{1520}(x-7)^{1729}$, cf.~\cite[Proposition 1]{Macaj-2010}. Since the coefficient of $x^{3245}$ in the polynomial $p_X(x)$ is $-116188800$ it follows that $X$ has 58094400 5-cycles. Given the fact that $X$ has degree 57 and order 3250, it is now straightforward to verify that $X$ also satisfies the equality.
\qed

Theorem~\ref{thm:main2} now follows immediately from Lemma~\ref{CL:charoddbound}.

\subsection{Even convex cycles}\label{SEC:even}

We next derive an upper bound for the number of even convex cycles, denoted with $\rho_e(G)$. The bound is similar to the above bound for $\rho_o(G)$.

It follows from the second part of Lemma \ref{LM:ccchar} that if $(v,v')$ is an even antipodal pair then $d_G(v,v') \geq 2$. Combining this with the fact that every even convex cycle $C$ yields $|C|/2$ antipodal pairs, gives the bound $$\rho_e(G) \leq \frac{ n(n-1)-2m}{g}\,.$$
While this bound is of the right order, it is not very sharp for sparse graphs. The next result establishes a better bound for graphs with a small cyclomatic number, that is, with a small $m-n+1$. 

\begin{lemma}
$\rho_e(G) \leq \bnd$. Moreover, equality holds if and only if $G$ is an even cycle.
\end{lemma}

\proof
We claim that every vertex $v \in V(G)$ lies in at most $m-n+1$ even antipodal pairs. Let $(v,v')$ be an antipodal pair of vertices from an even convex cycle $C$. Let $T$ be a BFS tree rooted at $v$. Lemma~\ref{LM:ccchar} implies that all the edges of $C$ are on $T$ with the exception of one edge $e$ that is incident with $v'$ on $C$. So for every vertex $v'$ that is antipodal with $v$ there is at least one edge $e$ not on $T$ that is adjacent to $v'$. This proves the claim. In total we therefore have at most $n(m-n+1)$ even antipodal pairs. In addition, every even convex cycle of length $2k$ yields $k$ antipodal pairs. Since we only need to count unordered pairs, we deduce that 
$$\rho_e(G) \leq \bnd\,.$$
For the equality part of the lemma, let $C$ be an even convex cycle of $G$. If $G = C$ then equality clearly holds. Otherwise, let $u$ be a vertex of $G$ that is not on $C$ and is adjacent to a vertex $v \in C.$ Let $v'$  be the antipodal vertex of $v$ on $C$. Then observe that $(u,v')$ is not an antipodal pair. Moreover, at least one edge that is incident with $v'$ on $C$ is not on a BFS tree rooted at $u$. We deduce that $u$ is contained in  less than $m-n+1$ even antipodal pairs which implies that $G$ has less than $\bnd$ even convex cycles.
\qed

\subsection{A combined inequality}

We finally combine the derived bounds for $\rho_o(G)$ and $\rho_e(G)$ into a single inequality for the number $\rho(G)$ of all convex cycles of $G$. The key insight is that graphs with the maximum number of convex cycles are homogeneous in the sense that they either contain only even or only odd convex cycles. The following lemma establishes this fact.

\begin{lemma}\label{LM:combbound}
$\rho(G) \leq \bnd$. Moreover, if $G$ contains an even convex cycle then the bound is attained if and only if $G = C_n$. 
\end{lemma}

\proof
Suppose that $C$ is an even convex cycle of $G.$ Let $v \in C$ and consider a BFS tree $T$ rooted at $v$. Let $v'$ be the antipodal vertex of $v$ with respect to $C.$ Let $e$ and $f$ be the edges of $C$ incident with $v'$. Then at least one of these two edges is not on $T$ and hence does not form an antipodal pair with $v$. This means that for every even convex cycle there is at least one less possible odd convex cycle which in turn implies that
$$\rho(G) \leq \bnd\,.$$ 

Suppose now that $G$ contains an even convex cycle $C$ and that $G \not= C_n$. Let $u \not \in C$ be a vertex of $G$ that is adjacent to a vertex $v \in C$. Let $v'$ be the antipodal vertex of $v$ on $C$ and consider a shortest $u,v'$-path $P_{uv'}$. We distinguish two cases and wish to show that the given configuration forbids the attainment of the bound. 

\medskip\noindent
{\bf Case 1.} $P_{uv'} \cap C \not= \emptyset$. \\
In this case $(u,v')$ is not an antipodal pair of an even convex cycle. Moreover at least one edge incident with $v'$ on $C$ is not in a BFS tree rooted at $u$ and also does not form an antipodal pair with $u.$ The latter fact implies that $\rho(G) < \bnd$.

\medskip\noindent
{\bf Case 2.}  $P_{uv'} \cap C = \emptyset$. \\
In this case the degree of $v'$ is at least $3$ and, because $C$ is convex, $|P_{uv'}| = d_G(v,v')$. It follows that a BFS tree $T$ rooted at $u$ does not contain the edges $e$ and $f$ that are on $C$ incident with $v$.  Moreover, none of these two edges forms an antipodal pair with $u$. Since $(u,v')$ is an antipodal pair of at most one even convex cycle, $u$ is contained in strictly less than $m-n+1$ antipodal pairs. Therefore the inequality for $\rho(G)$ is again not attained.  
\qed

Theorem~\ref{thm:main} now follows by combining Lemma~\ref{LM:combbound} with the results of Subsections~\ref{SEC:odd} and~\ref{SEC:even}.

\section{Concluding remark}

In this paper we have characterized the graphs in which the upper bound for the number of convex cycles $\bnd$ is attained. It turned out that there are not many such graphs. In might hence be interesting to study graphs that are close to this bound. A reasonable class of graphs in this respect could be generalized Moore graphs~\cite{Cerf-1974,Sampels-2004} as it appears that they contain many (even) convex cycles.

\section*{Acknowledgments}

This work has been financed by ARRS Slovenia under the grant
P1-0297 and within the EUROCORES Programme EUROGIGA/GReGAS of the
European Science Foundation. The second author is also with the Institute
of Mathematics, Physics and Mechanics, Ljubljana.


\begin{thebibliography}{99}

\bibitem{Bandelt-1996a}
  H.-J.~Bandelt, G.~Burosch, J.-M.~Laborde, 
  Cartesian products of trees and paths,
  J. Graph Theory 22 (1996) 347--356.

\bibitem{Bandelt-1996b}
  H.-J.~Bandelt, V.~Chepoi, 
  Graphs of acyclic cubical complexes,
  European J. Combin. 17 (1996) 113--120.

\bibitem{bannai-1973}
  E.~Bannai, T.~Ito, 
  On finite Moore graphs,
  J. Fac. Sci. Univ. Tokyo Sect. IA Math. 20 (1973) 191--208.

\bibitem{Biggs}
  N.~L.~Biggs, 
  Algebraic Graph Theory. Second Edition,
  Cambridge University Press, Cambridge, 1993. 

\bibitem{Bresar-2006}
  B.~Bre\v sar, S.~Klav\v zar, R.~\v Skrekovski, 
  Roots of cube polynomials of median graphs,
  J. Graph Theory 52 (2006) 37--50.

\bibitem{Cameron}
  P.~Cameron, 
  Permutation Groups, 
  Cambridge University Press, Cambridge, 1999.

\bibitem{Cerf-1974}
 V.~G.~Cerf, D.~D.~Cowan, R.~C.~Mullin, R.~Stanton, 
 Computer networks and generalized Moore graphs,
 in: Proceedings of the Third Manitoba Conference on Numerical Mathematics 
 (Winnipeg, 1973) (1974) 379--398.

\bibitem{Chepoi-1988}
  V.~Chepoi, 
  $d$-Convexity and isometric subgraphs of Hamming graphs,
  Cybernetics 1 (1988) 6--9.

\bibitem{Damerell-1973}
  R.~M.~Damerell, 
  On Moore graphs,
  Proc. Cambridge Philos. Soc. 74 (1973) 227--236.

\bibitem{Djokovic-1973}
  D.~Djokovi\'c, 
  Distance preserving subgraphs of hypercubes,
  J. Combin. Theory Ser. B 14 (1973) 263--267.

\bibitem{Egawa-1986}
  Y.~Egawa, 
  Characterization of the Cartesian product of complete graphs by convex subgraphs,
  Discrete Math. 58 (1986) 307--309.
  
\bibitem{Godsil}
  C.~Godsil, G.~Royle,
  Algebraic Graph Theory, 
  Springer, New York, 2001. 

\bibitem{Hammack-2011}
  R.~Hammack, W.~Imrich, S.~Klav\v zar, 
  Handbook of Product Graphs, Second Edition,
  CRC Press, Boca Raton, FL, 2011. 

\bibitem{Harary} 
  F.~Harary,  
  The determinant of the adjacency matrix of a graph, 
  SIAM Review 4 (1962) 202--210.

\bibitem{Hellmuth-2011}
  M.~Hellmuth, J.~Leydold, P.~F.~Stadler,
  Convex cycle bases and Cartesian products,
  manuscript, 2011.  

\bibitem{Hoffman-1960}
  A.~J.~Hoffman, R.~R.~Singleton, 
  On Moore graphs with diameters 2 and 3,
  IBM J. Res. Develop. 4 (1960) 497--504.

\bibitem{Klavzar-2012a}
  S.~Klav\v zar, S.~Shpectorov, 
  Convex excess in partial cubes,
  J. Graph Theory 69 (2012) 356--369.
  
\bibitem{Klavzar-2012b}
  S.~Klav\v zar, S.~Shpectorov, 
  Characterizing almost-median graphs II,
  Discrete Math. 312 (2012) 462--464.

\bibitem{Macaj-2010}
  M.~Ma\v caj, J. \v Sir\'{a}\v n,
  Search for properties of the missing Moore graph,
  Linear Algebra Appl. 432 (2010) 2381--2398.

\bibitem{Miller-2005}
  M. Miller, J. \v Sir\'{a}\v n,
  Moore graphs and beyond: A survey of the degree/diameter problem,
  Electron. J. Combin., Dynamic Survey DS14 (2005) 61 pp.

\bibitem{Polat-2009a}
  N.~Polat, 
  Netlike partial cubes. II. Retracts and netlike subgraphs,
  Discrete Math. 309 (2009) 1986--1998. 

\bibitem{Polat-2009b}
  N.~Polat, 
  Netlike partial cubes. IV. Fixed finite subgraph theorems,
  European J. Combin. 30 (2009) 1194--1204.

\bibitem{Rowlinson}
  P.~Rowlinson, I.~Sciriha,
  Some properties of the Hoffman-Singleton graph,
  Appl. Anal. Discrete Math. 1 (2007) 438--445. 

\bibitem{Sampels-2004}
  M.~Sampels, 
  Vertex-symmetric generalized Moore graphs,
  Discrete Appl. Math. 138 (2004) 195--202.
 
\bibitem{Siagiova-2012}
  J.~\v Siagiov\'{a}, J.~\v Sir\'{a}\v{n}, 
  Approaching the Moore bound for diameter two by Cayley graphs,
  J. Combin. Theory Ser. B 102 (2012) 470--473. 

\bibitem{Singleton} 
  R.~Singleton, 
  There is no irregular Moore graph, 
  Amer. Math. Monthly 75 (1968) 42--43.   

\bibitem{Vanden-1982}
  P.~Vanden Cruyce, 
  A characterization of the $n$-cube by convex subgraphs,
  Discrete Math. 41 (1982) 109--110. 
 
\bibitem{Vel-1993} 
  M.~L.~J.~van de Vel, 
  Theory of Convex Structures,
  North-Holland Publishing Co., Amsterdam, 1993.

\end{thebibliography}
\end{document}